\documentclass[10pt]{amsart}

\usepackage{latexsym,amsfonts,amssymb,epsfig,verbatim}
\usepackage{amsmath,amssymb,latexsym,graphics,textcomp,hyperref,enumerate} 
\usepackage[foot]{amsaddr}
\usepackage[utf8]{inputenc}
\usepackage{amssymb}
\usepackage{comment}
 \usepackage[T1]{fontenc}
 \usepackage{lmodern}
 \usepackage{graphicx} 
\usepackage{xspace} 
  
\usepackage[usenames,dvipsnames]{xcolor}

\usepackage{tikz}
\usepackage{pgfplots}

\usepackage{epigraph}
\usepackage{bm}
\usepackage{amsfonts}
\usepackage{amsmath}
\usepackage{enumerate}
\usepackage{eucal}
\usepackage{amsthm}
\usepackage{thmtools}
\usepackage[capitalize]{cleveref}
\newtheorem{theorem}{Theorem}[section]

\theoremstyle{definition}

\newcommand{\reals}{\mathbb{R}}

\newcommand{\mute}[1] {}

\def \grad{\operatorname{grad}}

\newcommand{\Int}{\operatorname{Int}}

\newcommand{\Fr}{\mathrm{Fr}}

\title[On the construction of geographical maps]{On the construction of geographical maps: Lagrange, Chebyshev, Darboux and Milnor}
\author{Hideki Miyachi}
\author{ Ken'ichi Ohshika}
\author{Athanase Papadopoulos}
\address{Hideki Miyachi,   
School of Mathematics and Physics,
College of Science and Engineering,
Kanazawa University,
Kakuma-machi, Kanazawa,
Ishikawa, 920-1192, Japan and Max-Planck-Institut für Mathematik, Vivatsgasse 7, 53111 Bonn, Germany}
\email{miyachi@se.kanazawa-u.ac.jp}  
\address{Ken'ichi Ohshika,
Department of Mathematics,
Gakushuin University,
Mejiro, Toshima-ku, Tokyo, Japan and Max-Planck-Institut für Mathematik, Vivatsgasse 7, 53111 Bonn, Germany}
  \email{ohshika@math.gakushuin.ac.jp}
\address{Athanase Papadopoulos,
Institut de Recherche Mathématique Avancée
(Université de Strasbourg et CNRS),
7 rue René Descartes,
67084 Strasbourg Cedex France and Max-Planck-Institut für Mathematik, Vivatsgasse 7, 53111 Bonn, Germany}
  \email{papadop@math.unistra.fr}
 
\begin{document}

\maketitle
%
%
    
 \begin{abstract}
Lagrange, Chebyshev, and Darboux, in 1779, 1856, and 1911, respectively,  wrote  articles all bearing the same title, \emph{On the Construction of Geographical Maps}. In 1969, Milnor wrote a paper in which he refers to Chebyshev's paper, of which he provides a new formulation and proof.
In the present article, we review the results of all these papers, explaining the main ideas they contain and pointing out connections between them. We give complete proofs of the statements by Darboux and Milnor, both of which aim to make explicit and provide a proof of Chebyshev's result, but whose contents are different. Although Chebyshev did not state explicitly what the word ``best" means, his conclusion, like that of Darboux and of Milnor, is that a best geographical map is characterised by the fact that its conformal factor is constant on the boundary of the region represented. Our statement and proof of Milnor's theorem work in a more general setting than the one he gives.

The final version of this paper will appear in the Handbook of Mathematics in the Arts and Sciences (second edition), ed. Bharath Sriraman, Springer, 2027.
  \bigskip


 \bigskip


\end{abstract}

\section{Introduction: The problem of mathematical geography}

A geographical map is a representation of a certain region of the surface of the Earth or of the Celestial globe (both generally considered to be a sphere, or, sometimes, for what regards the Earth, a spheroid) on a Euclidean plane.
The mathematical problem of constructing geographical maps is that of finding a map which satisfies certain predefined properties and has minimal distortion.
The variety of the available methods for constructing such maps depend on the precise properties required from the map, for example: sending meridians or parallels to line segments, or to arcs of circles (equicentric or not),  to arcs of ellipses, etc. Lagrange, in his work on the construction of geographical maps  \cite{Lagrange-memoir, Lagrange},  made the  remark that  from a most general point of view, ``one only has to draw the meridians
and parallels following a given arbitrary rule, and place the various locations with
respect to those lines as they are on the surface of the Earth with respect to the
circles of longitude and latitude''; see \cite[p.\ 313]{Lagrange-english} for an English translation of Lagrange's memoir. Lagrange solved in this memoir the problem in the case where the images of the meridians and parallels are circles or straight lines.

Besides this general fact about specifying the  images of meridians or parallels, a variety of other specific conditions may be desirable for a geographical map. For example, one may request that the map preserves area (up to a constant), or angles, or distances (up to scale) along some special parallel, or meridian, etc. Here, it is good to remember that it has been known since antiquity that it is impossible for a map from a domain on the sphere to the plane  to preserve distances up to scaling. A recent work on mathematical geography is the book \cite{Geography-Book}, in which all these questions are discussed in detail, and which contains English translations of works by eminent mathematicians on the problem of drawing geographical maps.
We also discuss the problem of drawing geographical maps in the article \cite{Handbook-article} contained the present Handbook.

The question of drawing geographical maps, as we have just addressed it (apart from the mention of the spheroid), dates back to Ptolemy, who wrote a treatise,  the \emph{Geography} \cite{Ptolemy-Geography}, on the representation of terrestrial maps, and another one,  the \emph{Mathematical composition}, better known under the name \emph{Almagest} \cite{Ptolemy-Almagest} (which, in fact, preceded  his work on terrestrial geography), on drawing celestial maps. Ptolemy's   \emph{Almagest} is familiar to mathematicians because it contains several results in plane and spherical geometry, including tables of spherical trigonometry (in the form of tables of chords), and the famous geometrical result known as Ptolemy's theorem.
It is interesting to highlight the fact that the same kind of questions that arise for drawing terrestrial maps  arise for celestial maps, making a connection between the Earth and the Heavens.

When Lagrange wrote his memoirs on geography,  two famous  conformal geographical projections were known. One is the stereographic projection, known since Greek antiquity,  and the other  is Mercator's projection.
Lambert in \cite{Lambert} gave a new kind of conformal projections, of the type called conical projections. He also considered a one-parameter family of projections of the sphere which are all conformal except possibly at the North pole, interpolating between the stereographic and the Mercator projections. Let us also note that there is a geographical map that goes by the name ``Lagrange projection", see \cite{Germain, Tissot-1879, Tissot-1880, Tsinger}. 
A certain number of geographical maps are compared from the point of view of their distortion in the two papers \cite{Handbook-article} and \cite{Geography-GB-OMP}.  

%
%
%

Let us note that all along the history of mathematics, problems in geography have motivated the study of geometrical problems.  We may quote Darboux, who is one of the main figures in the present article, from his address at the Rome 1908 ICM, titled
\emph{Les origines, les méthodes et les problèmes de la géométrie infinitésimale} (The
origins, methods and problems of infinitesimal geometry) 
\cite{ICM-1908}. In this address he declares:\footnote{The translations from the French are ours.} 
 
 \begin{quote}\small
Like many other branches of human knowledge, infinitesimal geometry
originated in the study of a problem posed by practical applications. The ancients had
already sought to obtain flat representations of the various regions of the
Earth, and they had adopted the very natural idea of projecting the surface of
our globe onto a plane. For a very long time, attention was focused exclusively on these methods
of projection, limiting oneself simply to studying the best ways to choose,
in each case, the viewpoint and the projection plane. It was a geometer of
the greatest insight, Lambert --- Lagrange’s highly esteemed colleague at the Academy of
Berlin ---, who, by pointing out for the first time a property common to
Mercator's maps, known as reduced maps, and to those produced by the stereographic projection,
was the first to consider the theory of geographical maps from a truly general perspective, and proposed, in its full scope, the problem of representing
the Earth’s surface on a plane while preserving the similarity of infinitely small elements.
This fascinating question, which gave rise to research by Lambert himself,
by Euler, and to two very important memoirs by Lagrange, was addressed for the first time by Gauss in all its generality.
 
\end{quote}

Darboux then goes on explaining how the problems of cartography led to many important developments in geometry made by Euler, Lagrange, Lambert, Gauss, Beltrami, Chebyshev and others. For the work of Chebyshev on geography, the reader may refer to
 \cite{Chebyshev-St-Petersburg1, Papadopoulos-Intelligencer, Tapia}.  Darboux himself was most interested in cartography. His textbook \cite{Darboux-Lecons} on the differential geometry of surfaces, which for several decades was the major reference on this topic in France, contains at several places references to the problem of drawing geographical maps. In particular, in \S 133 of his treatise  (p. 236-243 of Book II, Part I), Darboux solves the following problem inspired by  Lagrange's memoir on the construction of geographic maps. 

\medskip

\noindent {\bf Problem} (Lagrange). Considering the Earth as a spheroid of revolution, find all the geographic maps in which the meridians are represented by arcs of circles.

\medskip

To tackle this problem, Darboux first notes that on any surface of revolution, the system of meridians and parallels are two conjugate isothermal families, that is, that there exists a coordinate system $(u,v)$ in which the metric  is proportional to $du^2+dv^2$, and where the curves of the coordinate system satisfy the so-called Laplace--Darboux equation characterising conjugate nets.
   This implies that the images of the system of meridians and parallels by the map considered must also form two conjugate isothermal families.  But Darboux had already proved that (1) in the plane, a family of circles is isothermal only if the conjugate isothermal family is also a family of circles; (2) any isothermal system consisting of circles is obtained by a plane inversion, from one of two systems: (A) the system formed by curves which are coordinates of a system of rectangular rectilinear coordinates; (B) the system formed analogously using polar coordinates.
He then concludes the following:

\medskip

\emph{If one of the two families of meridians or parallels is represented by circles, then the same holds for the other family, and one obtains all the solutions of Lagrange's problem by assigning to the meridians or parallels either the lines of the system (A), or the coordinate curves of the system (B), up to  applying an arbitrary inversion of the plane.}

\medskip

Let us also note that in a memorial article read at the French Academy of Sciences on December 10, 1917, Picard writes the following, 
 regarding Darboux \cite{Picard-Darboux}: 
``[\ldots] Darboux had other plans as well. He wanted to write a book on the famous problem which gave rise to infinitesimal geometry --- the problem of geographical maps; he was captivated by both the elegance and the practical importance of this problem, and he had explored it in depth in his teaching. May his lecture notes help us reconstruct it!''

Lagrange, in his memoir on geography, refers to Lambert, writing: ``Lambert is the first who considered the theory of geographical maps
from the general point of view I just presented, and who consequently had the idea
of determining the lines of the meridians and of the parallels by the only condition
that all the angles made in the plane of the map are equal to the corresponding
angles on the surface of the globe.''
There is a well-known geographical map called the Lambert map. It belongs to the species of conformal conical map.
We have studied it in the paper \cite{Handbook-article}. Comparing it to other conical maps, we have shown that it is the best one among them, for what regards the drawing of the 18th-century Russian Empire, a subject in which Euler was heavily involved. 
Lagrange's work on geography is analysed in the paper \cite{Charitos}. 
Talking about mathematicians who were also geographers, we recall that Gauss, besides being a mathematician, was, officially, a geographer, and that several among his mathematical discoveries regarding the differential geometry of surfaces originate in problems in geography.
His paper \emph{General solution of the problem: to represent the parts of a given surface on another so that
the smallest parts of the representation shall be similar to the corresponding parts
of the surface represented} \cite{Gauss-Allgemeine}, published in 1825 and which contains important results on
conformal representations of simply-connected surfaces, was motivated by geographical
questions.

\section{Chebyshev's result and its formulation by Darboux and Milnor}
 
In 1856, Pafnuty Chebyshev gave two addresses, both titled \emph{On the construction of geographic maps} (the same title as Lagrange's two memoirs on geography), at the Imperial Academy of Saint Petersburg, see \cite{Chebyshev-Oeuvres, Chebyshev-Oeuvres-II, Chebyshev}.  In the mathematical  problem which he discusses in these memoirs, he assumes the Earth is spherical, but he notes that this is only for the purpose of simplifying the formulae, and that his method can easily be extended to all possible assumptions about the shape of the globe.
 From his work, we can extract the following result: 

\medskip

\noindent {\bf Claim} (Chebyshev).  \emph{Let $\Omega$ be a closed domain on the sphere, and consider the family of all conformal maps from  $\Omega$ to the plane $\reals^2$.
Then a conformal map is best among all these conformal maps if and only if its conformal factor is constant on $\Fr \Omega$.}
\medskip

 In fact, Chebyshev did not state  in precise terms the problem which he intended to solve, 
 but, as we shall recall later in this article, several mathematicians, including Darboux and Milnor, provided precise results and statements for slightly different theorems, attributing them to Chebyshev. For instance, the word ``best''  in the above statement was not defined precisely by Chebyshev, but it will be possible to  hypothesise its meaning, from the conclusion of the theorems, carefully stated by Darboux and by Milnor. Let us also note that Darboux, in his article \cite{Darboux-construction}, refers to the ``beautiful theorem of Chebyshev'', which he sets out to prove, but he states that this theorem ``is neither proved nor even clearly stated" by Chebyshev. 

 Milnor, in his article \cite{Milnor-cartography}, states and proves a theorem which he attributes to Chebyshev, which is yet a version different from Darboux’s. 

One of Chebyshev’s contributions was to note that Laplace’s equation $\Delta \phi = 0$  allows one to find, under suitable assumptions,  the most advantageous geographical representation, and this crucial observation was used later by Darboux and Milnor. 

 Chebyshev uses a notion of conformal factor, or ``scaling factor'', a notion which he attributes to Lagrange and for which he uses the French word ``rapport d'agran\-dissement'' (expansion ratio).
He gives the following formula for the conformal factor:
 \[ \lambda=\frac{\sqrt{f'(u+t\sqrt{-1}) F'(u-t\sqrt{-1}) }}
 {\displaystyle \frac{2}{e^u+e^{-u}}}
 .\]
  
To explain this notion, recall that two Riemannian metrics $ds^2$ and $d(s')^2$ on a surface $S$ are said to be \emph{conformal} if there exists a positive function $\lambda$ on $S$ such that $d(s')^2=\lambda^2 ds^2$. 
Given a conformal map from a surface $S$ (in the case considered by Chebyshev, $S$ is a subset of the sphere) to the plane, there is an associated positive function $\lambda$  such that if $d(s')^2$ is the pullback of the metric on the sphere by this map, we have $ds^2=\lambda^2 d(s')^2$. In this notation, $\lambda$ is the scaling factor. 

We are using the notion of Riemannian metric. Indeed, the proofs we give for the results contained in the  papers by Darboux and by Milnor are not restricted to metrics on a surface induced by an embedding in 3-space, but hold for arbitrary surfaces equipped with Riemannian metrics.

Chebyshev gave an answer to the problem of finding the best projection among the conformal ones, although his answer was not accompanied by sufficient mathematical justification. It is our aim in the rest of this paper to present the proofs of Chebyshev's result given by Darboux and Milnor respectively.

We shall state the results in a general form. We start by fixing some notation.

Instead of the sphere, we first consider a general surface $F$ equipped with a Riemannian metric.
Let $\Omega$ be a closed domain in $F$ and $f \colon \Omega \to \reals^2$ a conformal map from $\Omega$ to a subset of the plane.
Denote the Riemannian metric on $F$ by $ds^2$, and the pullback of the Euclidean metric under $f$ by $d(s')^2$.
Since $f$ is conformal, there is a function $\sigma_f$ on $\Omega$, which Darboux calls the infinitesimal scale function, characterised by the fact that for every $x \in \Omega$, we have \[d(s')^2(x)=\sigma_f(x)ds^2(x).\]
(This is equal to $\lambda^2$ in our previous notation.)
Darboux's interpretation of Chebyshev's claim is then the following.

\begin{theorem}[Darboux on Chebyshev]
\label{Darboux}
Let $f\colon \Omega \to \reals^2$ be a conformal map which minimises $\displaystyle \iint_{\Omega} ds^2(\grad(\log\sigma_g(x)), \grad(\log\sigma_g(x)) dA$ among all the conformal maps $g$ from $\Omega$ into $\reals^2$.
Then such a map $f$ exists and $\sigma_f(x)$ is constant on $\partial \Omega$.
\end{theorem}


Milnor's interpretation is as follows, restricted to the case when $F$ is a sphere.

\begin{theorem}[Milnor on Chebyshev]
\label{Milnor}
Let $\Omega$ be a closed domain on the sphere.
Let $f\colon \Omega \to \reals^2$ be a conformal map which minimises $\displaystyle\frac{\sup_{x\in \Omega} \sigma_g(x)}{\inf_{x\in \Omega} \sigma_g(x)}
$ among all conformal maps $g$ from $\Omega$ into $\reals^2$.
Then $\sigma_f(x)$ is constant on $\partial \Omega$, and such a map $f$ exists uniquely up to the similarity on $\reals^2$.
\end{theorem}

In \S4 of this paper, we shall prove a generalisation of Milnor's theorem.
We shall deal with a region on a more general surface, under the assumption that its Gaussian curvature is either non-negative or non-positive. The statement is the following.

\begin{theorem}
\label{general Milnor}
Let $\Omega$ be a closed domain on a surface with Riemannian metric whose Gaussian curvature is either non-negative or non-positive.
Let $f\colon \Omega \to \reals^2$ be a conformal map which minimises $\displaystyle\frac{\sup_{x\in \Omega} \sigma_g(x)}{\inf_{x\in \Omega} \sigma_g(x)}
$ among all conformal maps $g$ from $\Omega$ into $\reals^2$.
Then $\sigma_f(x)$ is constant on $\partial \Omega$, and such $f$ exists uniquely up to the similarity on $\reals^2$.
\end{theorem}

\section{Proof of Darboux's \cref{Darboux}}
\label{Th:Darboux}
Consider a general surface $F$ with a Riemannian metric $ds^2$, and a closed region $\Omega$ on $F$.
Express $ds^2$ using isothermal coordinates, so that we have
$$ds^2=\lambda^2 (du^2+dv^2).$$

 Note that Darboux used isothermal coordinates only in the second part of his argument.
As a matter of fact, in Darboux's time, the existence of isothermal coordinates were known only in the case of some special surfaces.

The gradient ${\rm grad}(\varphi)|_x\in T_x\Omega$ ($x\in \Omega$) of a $C^1$-function $\varphi$ on $\Omega$ is defined by
$$
d\varphi(X)=ds^2({\rm grad}(\varphi)|_x,X),
$$
 for all tangent vectors $X\in T_x\Omega$.
To be more concrete, this is expressed as
$$
\grad(\varphi)=\frac{1}{\lambda^2}
\left(
    \varphi_u\dfrac{\partial}{\partial u}+
    \varphi_v\dfrac{\partial}{\partial v}
\right).
$$
For a $C^2$-function $\varphi$ on $\Omega$, its Laplacian $\triangle_{ds^2}$ with respect to the metric $ds^2$ is expressed as
$$
\triangle_{ds^2}\varphi={\rm div}(\grad(\varphi))=
\frac{1}{\lambda^2}(\varphi_{uu}+\varphi_{vv}).
$$ 

Let $g(u,v)=(\xi,\eta)$ be a conformal map on $\Omega$.
Since
$$
g^*(du^2+dv^2)=\sigma_gds^2=(\sigma_g\lambda^2)(du^2+dv^2),
$$
the metric $(\sigma_g\lambda^2)(du^2+dv^2)$ on $\Omega$ is a flat metric, hence
\begin{equation}
\label{eq:equation_sigma_f}
0=\triangle_{ds^2}\log(\sigma_g\lambda^2)=\triangle_{ds^2}\log \sigma_g-K_{ds^2},
\end{equation}
where $K_{ds^2}$ is the Gaussian curvature of $ds^2$.

Suppose that the minimiser $f$ in \cref{Darboux} exists.
Let $g$ be another conformal map from $\Omega$ to $\mathbb{R}^2$.
Set $\psi=\log\sigma_g-\log \sigma_f$.
Then, from \eqref{eq:equation_sigma_f},
$$
\triangle_{ds^2}\psi=\triangle_{ds^2}(\log\sigma_g-\log \sigma_f)=-K_{ds^2}+K_{ds^2}=0.
$$
Hence, $\psi$ is harmonic on $\Omega$.

We set $\varphi=\log \sigma_f$.
Then we have  
\label{eq:I(g)}
\begin{align*}
I(g)&=\iint_\Omega ds^2(\grad(\log \sigma_g),\grad(\log \sigma_g))dA \\
&=\iint_\Omega \lambda^2\dfrac{(\log \sigma_g)_u^2+(\log \sigma_g)_v^2}{\lambda^4}\lambda^2dudv \\
&=\iint_\Omega((\log \sigma_g)_u^2+(\log \sigma_g)_v^2)dudv\\
&=\iint_\Omega((\log e^{\varphi+\psi})_u^2+(\log e^{\varphi+\psi})_v^2)dudv
\\
&=\iint_\Omega((\varphi+\psi)_u^2+(\varphi+\psi)_v^2)dudv
\\
&=I(f)+\delta I,
\end{align*}
 
where
$$
\delta I
=2\iint_\Omega (\psi_u\varphi_u+\psi_v\varphi_v)dudv
+\iint_\Omega (\psi_u^2+\psi_v^2)dudv.
$$
By Green's formula applied to the first term,
$$
\iint_\Omega (\psi_u\varphi_u+\psi_v\varphi_v)dudv
=-\int_{\partial \Omega} \varphi_g\frac{\partial \psi}{\partial n}ds,
$$
where $n$ is the inner normal vector to the boundary $\partial \Omega$.
Therefore, 
\begin{equation}
\label{eq:delta-I}
\delta I=-2 \int_{\partial \Omega} \varphi\frac{\partial \psi}{\partial n}ds+
\iint_\Omega (\psi_u^2+\psi_v^2)dudv.
\end{equation}

Since $I$ is assumed to take a minimum value at $f$,  the variation $\delta I$ is $\geq 0$.
We consider $\epsilon \psi$ instead of $\psi$ in \eqref{eq:delta-I}.  As $\epsilon \to 0$, we get the necessary condition that
\begin{equation}
\label{eq:Dalboux_condition}
\int_{\partial \Omega}\varphi \frac{\partial \psi}{\partial n}ds=0
\end{equation}
for all harmonic functions $\psi$ on $\Int \Omega$ that are differentiable on $\Omega$.

Conversely, suppose that $f$ satisfies \eqref{eq:Dalboux_condition} for all harmonic functions $\psi$ differentiable on $\Fr \Omega$. From \eqref{eq:delta-I}, we have
$$
\delta I=\epsilon^2\iint_\Omega (\psi_u^2+\psi_v^2)dudv>0
$$
for all harmonic functions $\psi$ on $\Int \Omega$ differentiable on $\Fr \Omega$.
Thus we see that Condition \eqref{eq:Dalboux_condition} is also a sufficient condition for $\delta I$ to be positive.

For a given smooth function $\eta$ on $\Fr \Omega$ with $\displaystyle\int_{\Fr \Omega}\eta ds=0$, we find a harmonic function $\psi$ which is smooth on $\Omega$ and satisfies the Neumann condition
$$
\frac{\partial \psi}{\partial n}=\eta
$$
(cf. \cite{Fo}).
Therefore, Condition \eqref{eq:Dalboux_condition} is equivalent to
$$
\int_{\Fr\Omega}\log \sigma_f\cdot \eta ds=0
$$
for all smooth function $g$ on $\Fr \Omega$ with $\int_{\Fr \Omega}\eta ds=0$.
By integrating by parts, this implies that $\log \sigma_f$ is constant on $\Fr \Omega$.
The converse trivially holds.
Thus, we have shown that the condition in the statement of  \cref{Darboux} is satisfied if and only if a conformal map $f$ minimises $I$.

Now, we turn to the existence of such a map $f$.
The problem is reduced to finding a function $\sigma$ satisfying
$$
\triangle_{ds^2}\log \sigma=K_{ds^2}
$$
on $\Int\Omega$ such that $\log \sigma$ is constant on $\partial \Omega$.
We set $\sigma =e^{-\mu}\lambda^2$.
Then $\sigma$ is a solution to the above problem if and only if the function $\mu$ is harmonic on the interior of $\Omega$, of class $C^1$ on $\Omega$ and satisfies
$\mu=2\log \lambda+c$ on $\partial \Omega$ for some $c\in \mathbb{R}$.
Such a function $\mu$ is obtained as the solution of the Dirichlet problem with the boundary value $2\log \lambda+c$ (see \cite{Fo}).

\section{Proof of Milnor's  \cref{general Milnor}}\label{Th:Milnor}
As we saw in the preceding section, Darboux minimises the quantity which we denote by $I(g)$. Instead, Milnor minimises a quantity  $\displaystyle\frac{\sup \sigma_g}{\inf \sigma_g}$, or, equivalently (by taking logarithms) $\sup \log\sigma_g -\inf \log \sigma_g$, which appears in (\ref{eq:milnor-comparison}) below.

Suppose first that the Gaussian curvature of the domain is non-negative. (We shall deal below with the case where it is non-positive.)
The proof in this case is the same as the one that Milnor gives in the case where the domain is a sphere.
We set $\phi_f(x)=\log \sigma_f(x)$ and $\phi_g(x)=\log\sigma_g(x)$.
It is well known that the differential equation $\triangle_{ds^2}\phi_f=K_{ds^2}$ has a unique solution satisfying the boundary condition $\phi_f(x)=0$ on $\partial \Omega$.
Let $g$ be another $C^2$-function satisfying $\triangle_{ds^2}\phi_g=K_{ds^2}$.
What we have to prove is  that
\begin{equation}
\label{eq:milnor-comparison}
\sup \phi_g-\inf \phi_g\ge \sup \phi_f-\inf \phi_f
\end{equation}
where equality holds if and only if
$$
\phi_g=\phi_f+\text{constant}.
$$

Since $\triangle_{ds^2}\phi_f\ge 0$, $\phi_f$ is subharmonic, and by the maximum principle, we have $\sup \phi_f=0$.
Set $c=\sup \phi_g$.
The difference $\phi_g-\phi_f$ satisfies
$\triangle_{ds^2}(\phi_g-\phi_f)=K_{ds^2}-K_{ds^2}=0$, which means that $\phi_g-\phi_f$ is harmonic.
Therefore, by the maximum principle again, there is a sequence $(x_i)_{i=1}^\infty$ in $\Omega$ converging to a point in $\Fr \Omega$ such that
$$
\lim_{i\to \infty}(\phi_g(x_i)-\phi_f(x_i))=\sup(\phi_g-\phi_f).
$$
Since $(x_i)$ tends to the boundary, we have  $\lim_{i\to \infty} \phi_f(x_i)=0$, and the left hand side is bounded by $c$.
This implies that  $\sup(\phi_g-\phi_f)\le c$, hence
\begin{equation}
    \label{eq:milner-comparison-2}
\phi_g(x)\le \phi_f(x)+c
\end{equation}
for all $x\in \Omega$. Thus we obtain
$$
\inf \phi_g\le \inf \phi_f+c.
$$
 Since $\sup \phi_f=0$, this is equivalent to \eqref{eq:milnor-comparison}.

Suppose that equality hold in \eqref{eq:milnor-comparison}.
Let $x_0\in \Omega$ be a point where $\phi_f$ takes its minimum.
From \eqref{eq:milner-comparison-2}, we have $\phi_g(x_0)\le \phi_f(x_0)+c$. If $\phi_g(x_0)<\phi_f(x_0)+c$, we have
$$
\inf(\phi_g)\le \phi_g(x_0)<\phi_f(x_0)+c=\inf(\phi_f)+\sup(\phi_g),
$$
which contradicts our assumption that equality holds.
Hence $\phi_g(x_0)=\phi_f(x_0)+c$, and $\phi_g-\phi_f$ achieves its maximum at $x=x_0\in\Int\Omega$. 
By the maximum principle, this means that  $\phi_g-\phi_f$ is constant.

We can deal with the case where $ds^2$ is non-positively curved using the same argument but exchanging $\sup$ and $\inf$.
Indeed we shall show that \eqref{eq:milnor-comparison} also holds under this assumption.

Let $f$ and $g$ be maps as above. 
Since $K_{ds^2}$ is non-positive in this case, $\phi_f$ and $\phi_g$ are super-harmonic whereas the difference $\phi_g-\phi_f$ is harmonic, as before.
This implies that $\inf \phi_f=0$. 
Now, set $c=\inf \phi_g$. Then, we have
$$
\inf(\phi_g-\phi_f)\ge c
$$
and $\phi_g(x)\ge \phi_f(x)+c$ for all $x\in \Omega$.
This means that $\sup(\phi_g)\ge \sup(\phi_f)+c$, which is equivalent to \eqref{eq:milnor-comparison}.
We can also deal with the case where the equality holds by the same argument as in the case of non-negative curvature.

Milnor, in his paper \cite{Milnor-cartography}, writes about Chebyshev's result: ``This result has been available for more than a hundred years, but to my
knowledge it has never been used by actual map makers." He does not mention Darboux' paper in which the latter gave his own version of Chebyshev's result. It seems likely that  when he wrote his paper \cite{Milnor-cartography}, Milnor was not aware of Darboux' paper.

\section{In guise of a conclusion}
Let us end this paper with some words of conclusion. Three ideas  emerge from this work. The first one is the role of mathematicians in the art of drawing geographical maps. The second  is that mathematics is a collective endeavour: several  people, over the course of several centuries, have raised, considered, reflected on, explained, and refined the same question, each contributing his own ideas and perspectives. The third  is that mathematics is timeless. The ideas, language, and techniques which our predecessors developed centuries ago are still the ones which interest us and which we use today.

\bigskip

\noindent{\bf Acknowledgement.} The authors are grateful to Bharath Sriraman who read an early version of this paper and made useful comments.

\bibliography{Tch}
\bibliographystyle{plain}
\end{document}